\documentclass[11pt]{article}
\usepackage[T1]{fontenc}
\usepackage{lmodern}
\usepackage{amsmath,amssymb,amsthm}
\usepackage[margin=1in]{geometry}
\usepackage{booktabs}
\usepackage[hidelinks]{hyperref}
\usepackage{microtype}
\newtheorem{theorem}{Theorem}[section]
\newtheorem{lemma}[theorem]{Lemma}
\newtheorem{proposition}[theorem]{Proposition}
\newtheorem{corollary}[theorem]{Corollary}
\theoremstyle{definition}
\newtheorem{definition}[theorem]{Definition}
\theoremstyle{remark}

\newcommand{\ee}{\mathbf 1}
\newcommand{\RR}{\mathbb R}
\newcommand{\MIN}{\mathsf{MIN}}
\newcommand{\supp}{\operatorname{supp}}
\newcommand{\argmin}{\operatorname*{arg\,min}}
\newcommand{\argmax}{\operatorname*{arg\,max}}
\newcommand{\diag}{\operatorname{diag}}

\newcommand{\T}{^{\mathsf T}}
\title{Minimum-coordinate elimination on a simplex:\newline
polynomial first-step safety and coNP-complete correctness}
\author{Arseny Kosachev\\[0.4em]
\small Department of Operations Research\\
\small Faculty of Computational Mathematics and Cybernetics\\
\small Lomonosov Moscow State University\\
\small Moscow, Russia}
\date{October 23, 2026}

\begin{document}
\maketitle
\begin{abstract}
We study a deletion rule for strictly convex quadratic minimization on a
simplex: compute the affine minimizer on the current face, permanently
remove a most-negative coordinate, and repeat. Correctness is required
for every linear term and every resolution of ties. We characterize
first-step safety by a finite maximum principle for the constrained
inverse, and correctness from every initial face by principal-inverse
row inequalities. Correctness from the full simplex is a different
property. The two properties coincide with at most four coordinates,
but an explicit rational five-coordinate example separates them. An
exact reachability condition explains the separation: boundary
equalities in the constrained inverse restrict the tasks that can be
entered on a facet. Away from these equalities, full-start and hereditary
correctness coincide. We also determine the recognition complexity.
First-step safety is decidable in polynomial bit time, by combining
classical strict $n$-step-vector recognition with diagonal regularization.
In contrast, full-start and hereditary correctness are coNP-complete,
even on a reduction family with strict column diagonal maxima on every
face. The proofs include the non-strict boundary and all affine anchors;
explicit rational certificates for the finite examples are given in the text.
\end{abstract}

\section{Introduction}
An affine minimizer of a strictly convex quadratic on a simplex is easy
to compute, but it may have negative coordinates. A natural elimination
rule removes a most-negative coordinate and solves the smaller affine
problem. This rule terminates after finitely many linear solves. Its
termination at a feasible point, however, does not by itself establish
optimality: an earlier deletion may have removed a coordinate needed by
the constrained optimum. The question studied here is which fixed
Hessians make every such deletion safe for every linear term.

\subsection{Problem and deletion rule}
Let $P=\{1,\ldots,n\}$ and let $A\in\RR^{n\times n}$ be symmetric positive
definite. For $b\in\RR^n$, consider
\begin{equation}\label{eq:qp}
 \min_{x\in\Delta_P} f_b(x),\qquad
 f_b(x)=\tfrac12x\T A x-b\T x,\qquad
 \Delta_P=\{x\ge0:\ee\T x=1\}.
\end{equation}
For a nonempty $J\subseteq P$, write $\Delta_J=\{x\in\Delta_P:
x_{P\setminus J}=0\}$. Its affine hull is defined by the same equalities,
without the inequalities $x_J\ge0$. The affine minimizer is unique.

The algorithm $\MIN$ starts from a specified nonempty $J$. It computes
the affine minimizer $y^J$ on that face. If $y^J_J\ge0$, it returns $y^J$.
Otherwise it chooses any
\[
 h\in\argmin_{i\in J}y^J_i
\]
and replaces $J$ by $J\setminus\{h\}$. Thus only strictly negative
coordinates are removed, removed coordinates never return, and ties may
be resolved arbitrarily. The algorithm terminates after at most $|J|-1$
deletions, since its affine coordinates sum to one.

\begin{definition}
A matrix $A$ has \emph{full-start correctness} if, for every $b$ and every
tie resolution, $\MIN$ starting from $P$ returns the minimizer of
\eqref{eq:qp}. It has \emph{hereditary correctness} if the analogous
statement holds from every nonempty initial face $J$, with the target
being the minimizer over $\Delta_J$.
\end{definition}

The target in the second definition is the optimum on the initial face;
an arbitrary face need not contain the optimum on $\Delta_P$.

\subsection{Results and organization}
The distinction between these two quantifiers is substantive.
\begin{theorem}\label{th:main}
Full-start correctness and hereditary correctness coincide for $n\le4$.
For $n=5$, there is a rational positive definite matrix with full-start
correctness but without hereditary correctness. All these statements
include every resolution of ties.
\end{theorem}

The matrix and its complete rational certificate are given in
Section~\ref{sec:example}. The obstruction is a boundary phenomenon in
the constrained inverse, rather than a numerical ambiguity. A generic
equivalence is proved in Section~\ref{sec:small}.

The remaining results identify the relevant matrix conditions and their
recognition complexity:
\begin{enumerate}
\item A finite maximum principle characterizes universal safety of the
first deletion (Theorem~\ref{th:first}). Applying it to every face gives
an exact criterion in terms of the negative off-diagonal parts of
principal inverses (Theorem~\ref{th:hereditary}).
\item A reverse-insertion formula characterizes exactly which affine
tasks can follow the first deletion (Theorem~\ref{th:guarded}). It yields
the five-coordinate separation and the generic equivalence
(Theorem~\ref{th:generic}).
\item First-step safety is recognizable in polynomial bit time, including
zero principal-inverse components (Theorem~\ref{th:fmp-poly}). Full-start
and hereditary correctness are both coNP-complete to recognize
(Theorems~\ref{th:hereditary-hardness} and~\ref{th:fullstart-hardness}).
The same rational reduction proves both hardness statements and controls
all choices of affine anchor.
\end{enumerate}
Sections~\ref{sec:constrained}--\ref{sec:reachable} establish the matrix
criteria. Sections~\ref{sec:example}--\ref{sec:small} prove the separation,
its dimensional minimality, and the generic contrast.
Sections~\ref{sec:recognition}--\ref{sec:hereditary-hardness} address
recognition. Appendix~\ref{app:basic-failure} records a small example of
failure of the unrestricted deletion rule, separately from the
full-start versus hereditary distinction. Appendix~\ref{app:return}
derives a cycle of an unsafeguarded return modification from an old
complementarity example and gives a small-dimensional termination bound.
Cycling is not part of the novelty claim for this article.

\subsection{Relation to projection and earlier procedures}
If the columns of an invertible matrix $V$ are simplex vertices and $p$
is a query point, then projection minimizes
$\tfrac12\|Vx-p\|^2$ on $\Delta_P$. This has $A=V\T V$ and
$b=V\T p$; varying $p$ permits every $b$. Conversely every positive
definite $A$ admits such a realization. The question is therefore also
one about projecting onto an arbitrary simplex.

The most-negative-variable selection rule has earlier precedents.
Clark's restricted least-squares method \cite[Sections~2.1 and 2.5]{clark}
repeatedly fixes a most negative free variable at zero. On reaching a
feasible point, that method checks full optimality conditions and may
return variables to the free set. Its branch-existence result does not
say that every most-negative choice can be made irreversible. The
algorithm studied here is precisely the deletion-only rule defined
above; neither that selection rule nor active-set reintroduction is
being proposed as a new algorithm.

Cevikalp, Triggs, Jurie, and Polikar describe the repeated elimination
of the most negative affine weight in their convex-hull procedure
\cite[Section~2.2]{cevikalp}. Interpreting that passage as universal
correctness of the literal deletion-only rule is disproved by the exact
three-vertex example in Appendix~\ref{app:basic-failure}. This observation
concerns that elimination rule, rather than the other components or
empirical results of \cite{cevikalp}.
Theorem~\ref{th:main} addresses a different issue: whether a matrix on
which the rule is universally correct must retain that property on faces.

The maximum principles and principal inverse conditions used below are
related to established potential-theoretic and complementarity concepts.
Frostman-type maximum principles occur in finite-network potential theory
\cite{bendito}. Nonnegative principal solutions with a fixed right-hand
side are the non-strict form of conditions appearing in the theory of
$n$-step vectors \cite{pang,adler}. Those conditions should not be confused
with correctness of the present deletion rule: complementary pivot
algorithms have their own pivot tests and transformed covering vectors.
The row inequality in Theorem~\ref{th:hereditary} also occurs as a matrix
condition in \cite{kager}. We give all proofs needed here and make no
priority claim for these general matrix ingredients or for the elementary
projection counterexample in the appendix.

Structured sufficient conditions address different quantifiers and
different deletion rules. For the homogeneous objective $b=0$, absolute
row diagonal dominance of an SPD $A$ suffices for most-negative deletion.
The relevant enclosure and certified batch bounds follow from the
interval complementarity theorem of
Alefeld--Wang--Shen~\cite[Corollary~2.2]{alefeld}; they are not claimed as
new results here. Simultaneously deleting every negative coordinate is
a stronger rule. For example, when $A^{-1}$ is a symmetric nonsingular
M-matrix, the homogeneous simplex problem reduces to the established
block-elimination procedure for an M-matrix complementarity problem
\cite[Theorems~2--3 and Example~3]{fang}. These special classes do not
replace the universal-in-$b$ criteria for the specific $\MIN$ rule
studied in this article.

The parametric Hessian algorithms of Best--Caron and Best--Chakravarti
also follow optimal solutions with both primal and dual event tests.
In the complete 1984 Best--Caron report \cite[Algorithm~D and
Theorem~1]{bestcaron}, the initial problem is already solved and active
constraints may be added or released; its convergence theorem assumes
strict initial multipliers and unambiguous events. Theorem~3.1 of
Best--Chakravarti \cite{bestchak} gives an $O(n^2)$ arithmetic bound for
their specified continuation method on a diagonal-plus-rank-one model.
Neither statement addresses all linear terms and all tie resolutions of
the permanent minimum-coordinate rule for a fixed arbitrary SPD matrix.
The Best--Caron theorem locators refer to their 1984 report.

Residualized Gram matrices also occur in the partial positive-cone
assumption for sequential Lasso \cite[Section~3.1, Assumption~A2]{luochen}.
That assumption uses signal-dependent selected sets and supports a
statistical consistency result. With fixed column orientations, our
hereditary criterion instead requires the non-strict inverse row-sum
inequality for every pair of disjoint eliminated and retained index sets,
whose union is proper for the centered $Q$; this follows from
\eqref{eq:subinverse}. We do not identify these different quantifiers.
The signed Lasso positive-cone condition includes every diagonal sign
choice; its recognition reduces to strict absolute diagonal dominance
of the full inverse Gram matrix \cite[Theorem~3]{duan}. In contrast,
\eqref{eq:eta} discards positive off-diagonal entries, and the stronger
signed condition is not an equivalent recognition test here.
Another polynomial recognition result is Morris's test for hidden
positive row diagonally dominant matrices \cite{morrishprdd}: it seeks
factors $C,B$ with $AC=B$ and strictly absolutely dominant rows in both
factors. That existential factorization problem is distinct from the
universal principal-inverse row conditions imposed here. We do not
transfer its recognition algorithm to those conditions.

\subsection{Scope of the contribution and recognition complexity}
The polynomial recognition ingredient for strict $n$-step vectors is
classical: Morris \cite[Theorem~6.1]{morris} states a linear-feasibility
characterization, crediting Pang and Morris--Lawrence for its two
directions. Section~\ref{sec:recognition} supplies an explicit rational
diagonal regularization to handle zero components for the SPD and
centered matrices in this problem. We present that result as a deduction
from classical recognition theory, not as a new general polynomial
algorithm for complementarity or a priority claim for weak recognition.
The publisher's abstract of Pang--Chandrasekaran~\cite{pangchandra}
already announces a polynomial testing procedure. Non-strict
$n$-step inequalities already
occur in Pang's report \cite{pang}, in the account of Chu
\cite[p.~463]{chu}, and in the bounded-variable analysis of
Pang--Han \cite[equation~(2.3) and Algorithm~I]{panghan}. Those pivot
algorithms assume a suitable covering vector and use parameter or ratio
tests; their iteration bounds do not certify arbitrary irreversible
minimum-coordinate choices.

The claims developed here concern instead the exact quantifiers of
$\MIN$: safety for all linear terms, correctness after every permitted
tie, and the distinction between full initialization and every initial
face. The boundary reachability restriction and its five-coordinate
separation are not conclusions of the cited $n$-step or positive-cone
theorems. The hardness statements classify recognition of these
particular correctness properties. They do not assert hardness of
solving the underlying convex quadratic program, which is polynomially
solvable in the rational bit model. Nor does classical hardness of
recognizing a general P-matrix prove them: the input matrices here are
SPD, and the reduction must protect every other row and every affine
anchor. These protections are proved explicitly in
Section~\ref{sec:hereditary-hardness}.

Two companion manuscripts by the same author address different
questions. The local analysis \cite{kosachevlocal} concerns correctness
of arbitrary negative-coordinate deletion throughout a neighborhood of
a fixed KKT point, including degeneracy. The Schur-complement study
\cite{kosachevschur} concerns conditional off-diagonal signs and selection
of a set whose elimination protects the remaining block. The present
hardness proof shares signed-graph and residual-Gram ideas with that
study; it gives its own reduction because principal-inverse row margins
and all affine anchors require additional estimates. All the results
used to prove this paper's main theorems are established here; neither
companion manuscript is an external assumption. Both are unpublished
manuscripts, and no arXiv identifiers are assigned to them here.

\section{Constrained inverses and a maximum principle}\label{sec:constrained}
For nonempty $J$, let $A_J=A_{JJ}$ and put
\begin{equation}\label{eq:qdef}
 Q_J=A_J^{-1}
 -\frac{A_J^{-1}\ee_J\ee_J\T A_J^{-1}}
 {\ee_J\T A_J^{-1}\ee_J}.
\end{equation}
This is the derivative of $y^J_J$ with respect to $b_J$.
It is positive semidefinite with kernel $\operatorname{span}(\ee_J)$.
For $|J|\ge2$, every nonempty proper principal submatrix is positive
definite: a vector supported on a proper subset cannot be a nonzero
multiple of $\ee_J$.
Write $Q=Q_P$.

We use the residual convention
\begin{equation}\label{eq:kkt}
 r=Ax^*-b-\nu\ee\ge0,\qquad x^*\ge0,\qquad
 \ee\T x^*=1,\qquad x_i^*r_i=0.
\end{equation}
These conditions characterize the unique optimum. If that optimum is
supported on $J$, its affine minimizer satisfies
\begin{equation}\label{eq:affine}
 y^J_J=x^*_J-Q_Jr_J.
\end{equation}
For example, $Q_JA_Jd=d$ whenever $\ee_J\T d=0$; subtracting affine
stationarity from \eqref{eq:kkt} proves \eqref{eq:affine}.

\begin{lemma}\label{lem:schur}
If $D\subsetneq J$ is nonempty and $I=J\setminus D$, then
\begin{equation}\label{eq:schur}
 Q_I=(Q_J)_{II}-(Q_J)_{ID}(Q_J)_{DD}^{-1}(Q_J)_{DI}.
\end{equation}
In particular, if $D=\{h\}$, $q=(Q_J)_{hh}$, and
$c_i=(Q_J)_{ih}/q$, then
\begin{equation}\label{eq:delete}
 y^I_i=y^J_i-c_i y^J_h\qquad(i\in I).
\end{equation}
\end{lemma}
\begin{proof}
Vary only $b_D$ until the affine solution has coordinates $y_D=0$.
The change in those coordinates is $(Q_J)_{DD}\,\delta b_D$, so the
required change is $\delta b_D=-(Q_J)_{DD}^{-1}y_D$.
The resulting solution is stationary on $I$ for the unchanged $b_I$;
it is therefore the affine minimizer on $I$. Its coordinates are
$y_I-(Q_J)_{ID}(Q_J)_{DD}^{-1}y_D$. Differentiating this identity with
respect to $b_I$ gives \eqref{eq:schur}. The singleton case gives
\eqref{eq:delete}.
\end{proof}

For a symmetric matrix $M$, define the following finite maximum principle:
\begin{equation}\label{eq:fmp}
 \max_i(Mr)_i=\max_{i\in\supp r}(Mr)_i
 \qquad\text{for every }r\ge0,\ r\ne0.
\end{equation}
No assertion about signed loads is part of this definition.

\begin{theorem}\label{th:first}
For fixed $A$, the first step of $\MIN$ is safe for every $b$ and every
tie resolution if and only if $Q$ satisfies \eqref{eq:fmp}.
Here safe means that the removed coordinate is zero at the optimum.
\end{theorem}
\begin{proof}
Suppose FMP holds and $x_i^*>0$, $y_i<0$. Then $r_i=0$ and
$(Qr)_i=x_i^*-y_i>0$. Some $j\in\supp r$ satisfies
$(Qr)_j\ge(Qr)_i$, so
\[
 y_j=-(Qr)_j\le-(Qr)_i< x_i^*-(Qr)_i=y_i.
\]
Thus $i$ cannot be a minimizing coordinate, even in a tie. If $r=0$,
the affine optimum is already feasible.

Conversely, let $r\ge0$ violate FMP. Set $O=\argmax_i(Qr)_i$ and
$\beta=\max_i(Qr)_i$. Then $O\cap\supp r=\varnothing$.
The support of $r$ is proper, and positive definiteness of its principal
block gives $r\T Qr>0$, hence $\beta>0$. Choose $x^*$ uniform on $O$
and set $b=Ax^*-tr$. The pair $(x^*,tr)$ satisfies \eqref{eq:kkt}
with $\nu=0$. For sufficiently large $t$, every minimum coordinate of
$y=x^*-tQr$ lies in $O$ and is negative. Thus every first choice removes
a positive optimal coordinate.
\end{proof}

\begin{theorem}\label{th:inverse}
For a symmetric positive definite matrix $M$, FMP is equivalent to
\begin{equation}\label{eq:invspd}
 M_{TT}^{-1}\ee_T\ge0\quad\text{for every nonempty }T,
\end{equation}
including the full index set. For $Q\succeq0$ with
$\ker Q=\operatorname{span}(\ee)$, it is equivalent to
\begin{equation}\label{eq:invcenter}
 Q_{TT}^{-1}\ee_T\ge0\quad\text{for every nonempty proper }T.
\end{equation}
\end{theorem}
\begin{proof}
For necessity, choose a set $T$ of minimum cardinality with a negative
component $i$ of $M_{TT}^{-1}\ee_T$. A singleton cannot be such a set.
Put $R=T\setminus\{i\}$ and $a=M_{RR}^{-1}\ee_R\ge0$. Block inversion
gives
\[
 (M_{TT}^{-1}\ee_T)_i
 =\frac{1-M_{iR}a}{M_{ii}-M_{iR}M_{RR}^{-1}M_{Ri}}.
\]
The denominator is positive, so $M_{iR}a>1$, while
$M_{RR}a=\ee_R$. Extending $a$ by zero violates FMP. The same argument
applies to proper principal blocks of $Q$.

For sufficiency in the positive definite case, use induction on the
matrix size. Proper principal blocks have FMP by induction. A load of
full support requires no proof. For a nonzero load supported on a proper
set $R$, the maximum potential on $R$ is positive, by positive energy;
normalize that maximum to one. Consider
\[
 \mathcal P_R=\{v\ge0:M_{RR}v\le\ee_R\}.
\]
This polytope is compact, since
$\lambda_{\min}(M_{RR})\|v\|^2\le v\T M_{RR}v\le\ee_R\T v$.
At a nonzero vertex $v$ with actual support $T\subseteq R$, all
inequalities on $T$ must be equalities:
\begin{equation}\label{eq:vertex}
 M_{TT}v_T=\ee_T.
\end{equation}
Indeed, FMP for $M_{RR}$ makes the inequalities on $R\setminus T$
redundant for nonnegative vectors supported on $T$. If fewer than $|T|$
of the remaining inequalities are active, a nonzero direction preserving
all active equalities permits sufficiently small perturbations of both
signs, with positive coordinates on $T$. This contradicts extremality.

For any $i\notin R$, applying \eqref{eq:invspd} and the block inverse
formula to $T\cup\{i\}$ gives $M_{iT}v_T\le1$. The zero vertex also
satisfies this inequality. By convexity it holds throughout $\mathcal P_R$,
which proves FMP for $M$.

In the centered case, all proper blocks $Q_{RR}$ satisfy the positive
definite result just proved. The same polytope argument works unless
$T\cup\{i\}$ is the full index set. In that case $Q\ee=0$ yields
$Q_{iT}=-\ee_T\T Q_{TT}$, hence
$Q_{iT}Q_{TT}^{-1}\ee_T=-|T|\le1$. This supplies the remaining
inequality. Zero components in \eqref{eq:invspd} or
\eqref{eq:invcenter} cause no difficulty.
\end{proof}

\section{Correctness on all faces and reachable facets}\label{sec:reachable}
For a nonempty proper $T\subset P$, write $H_T=Q_{TT}^{-1}$ and define
\begin{equation}\label{eq:eta}
 \eta_i(T)=(H_T)_{ii}
 +\sum_{\substack{j\in T\setminus\{i\}\\(H_T)_{ij}<0}}(H_T)_{ij}.
\end{equation}

\begin{theorem}\label{th:hereditary}
The following are equivalent:
\begin{enumerate}
\item $A$ has hereditary correctness;
\item $Q_J$ satisfies FMP for every nonempty face $J$;
\item $\eta_i(T)\ge0$ for every nonempty proper $T\subset P$ and $i\in T$.
\end{enumerate}
\end{theorem}
\begin{proof}
The first two statements are equivalent by Theorem~\ref{th:first} and
induction along the deletions. Safe deletions preserve the initial
optimum; termination at a feasible affine minimizer then returns it.

Let $D=P\setminus J$ and $\varnothing\ne R\subsetneq J$. The block
inverse formula and Lemma~\ref{lem:schur} give
\begin{equation}\label{eq:subinverse}
 \bigl((Q_J)_{RR}\bigr)^{-1}
 =\bigl(Q_{R\cup D,R\cup D}^{-1}\bigr)_{RR}.
\end{equation}
Every inverted block is proper and positive definite. By
Theorem~\ref{th:inverse}, assertion 2 is therefore equivalent to
\[
 \sum_{j\in R}(H_T)_{ij}\ge0
 \qquad(i\in R\subseteq T\subsetneq P).
\]
For fixed $T$ and $i$, the smallest such row sum includes precisely
all negative off-diagonal entries and omits all positive ones. That
smallest value is \eqref{eq:eta}.
\end{proof}

We next describe exactly which tasks on a facet can follow a first
deletion. This is the distinction needed for full-start correctness.

\begin{lemma}[Reverse insertion]\label{lem:insertion}
Let $I=P\setminus\{h\}$, $q=Q_{hh}>0$, and $c_i=Q_{ih}/q$.
Fix $b_I$ and let $z$ be its affine minimizer on $I$. For any $s\in\RR$
there is a unique $b_h(s)$ for which the full affine minimizer satisfies
\begin{equation}\label{eq:insertion}
 y_h=s,\qquad y_i=z_i+c_i s\ (i\in I),\qquad
 b_h(s)=b_h(0)+s/q.
\end{equation}
If $s\to-\infty$, the optimum on $I$, extended by a zero at $h$,
is also the full optimum for all sufficiently negative $s$.
\end{lemma}
\begin{proof}
The derivative of $y_h$ with respect to $b_h$ is $q>0$. At $y_h=0$,
the other coordinates equal $z$, so Lemma~\ref{lem:schur} or linearity
gives \eqref{eq:insertion}. If $x^I,\nu_I$ solve the constrained problem
on $I$, their new residual is
$r_h=A_{hI}x^I-b_h(s)-\nu_I\to+\infty$. All other KKT conditions
remain unchanged, proving the last assertion.
\end{proof}

Suppose now that $Q$ has FMP. Applying it to the coordinate load $e_h$
gives $c_i\le1$. Define the equality set
\begin{equation}\label{eq:equalityset}
 C_h=\{i\in P\setminus\{h\}:Q_{ih}=Q_{hh}\}.
\end{equation}

\begin{lemma}[Exact reachability condition]\label{lem:reach}
Assume FMP for $Q$ and fix $b_I$ on $I=P\setminus\{h\}$, with affine
minimizer $z$. There exists an extension $b_h$ that permits $h$ as a first
$\MIN$ deletion if and only if
\begin{equation}\label{eq:guard}
 z_i\ge0\quad(i\in C_h).
\end{equation}
When this condition holds, the extension can also make $h$ strictly
inactive at the true optimum while preserving the optimum on $I$.
\end{lemma}
\begin{proof}
For a first deletion of $h$ we need $s<0$ and
$y_i-y_h=z_i+s(c_i-1)\ge0$. For $i\in C_h$ this is exactly
\eqref{eq:guard}, proving necessity. For other indices $c_i<1$, so all
remaining inequalities are strict for sufficiently negative $s$.
The last assertion follows from Lemma~\ref{lem:insertion}. If some
$z_i=0$ on $C_h$, a tie is permitted by the definition of $\MIN$.
\end{proof}

\begin{theorem}[Characterization by reachable facets]\label{th:guarded}
For $n\ge2$, $A$ has full-start correctness if and only if both of the
following hold:
\begin{enumerate}
\item $Q$ satisfies FMP;
\item for every $h$, the algorithm starting on $I=P\setminus\{h\}$
is correct for every $b_I$ whose initial affine minimizer satisfies
$z_{C_h}\ge0$, and for every tie resolution.
\end{enumerate}
\end{theorem}
\begin{proof}
Full-start correctness implies assertion 1 by Theorem~\ref{th:first}.
If assertion 2 fails, Lemma~\ref{lem:reach} extends the failing task so
that $h$ is a permitted first deletion and the facet optimum remains the
full optimum. The same failing continuation then contradicts full-start
correctness.

Conversely, an initially feasible affine minimum is already optimal.
Otherwise, any first deletion $h$ is safe by assertion 1, so its facet
contains the full optimum. The affine minimizer $z$ on that facet obeys
$z_i=y_i-c_i y_h=y_i-y_h\ge0$ for $i\in C_h$.
Assertion 2 applies to every continuation and returns the facet optimum,
which equals the full optimum.
\end{proof}

\section{A rational counterexample with five coordinates}\label{sec:example}
Define
\begin{equation}\label{eq:ts}
 T=\frac1{23}\begin{pmatrix}32&12&-44\\12&16&-28\\-44&-28&72\end{pmatrix},
 \qquad
 S=\begin{pmatrix}
 32/23&12/23&-44/23&0\\
 12/23&39/23&-28/23&-1\\
 -44/23&-28/23&72/23&0\\
 0&-1&0&1
 \end{pmatrix}.
\end{equation}
Thus $S=\diag(T,0)+(e_2-e_4)(e_2-e_4)\T$. Put
\begin{equation}\label{eq:qconstruction}
 c=(-9/23,-12/23,-25/23,1)\T,\qquad q=1/2,
 \qquad Q=\begin{pmatrix}S+qcc\T&qc\\qc\T&q\end{pmatrix}.
\end{equation}
Equivalently,
\begin{equation}\label{eq:Qexplicit}
 Q=\frac1{1058}\begin{pmatrix}
 1553&660&-1799&-207&-207\\
 660&1938&-988&-1334&-276\\
 -1799&-988&3937&-575&-575\\
 -207&-1334&-575&1587&529\\
 -207&-276&-575&529&529
 \end{pmatrix}.
\end{equation}
The upper $2\times2$ block of $T$ is positive definite, with determinant
$16/23$, and $T\ee=0$. Hence $T\succeq0$ with kernel
$\operatorname{span}(\ee_3)$. The representation of $S$ gives the
corresponding facts for $S$ and $\ee_4$. Since $\ee\T c=-1$, the quadratic
form of $Q$ is
\[
 x\T Sx+q(c\T x+x_5)^2.
\]
Consequently $Q\succeq0$ and $\ker Q=\operatorname{span}(\ee_5)$.
Every nonempty proper principal block is positive definite.

The associated positive definite matrix is
\begin{equation}\label{eq:Aexplicit}
 \begin{split}
 A&=(Q+\ee\ee\T/5)^{-1}\\
 &=\frac1{2300}\begin{pmatrix}
 3026&-1804&841&-1804&2041\\
 -1804&3716&36&3716&-3364\\
 841&36&956&36&431\\
 -1804&3716&36&6016&-5664\\
 2041&-3364&431&-5664&8856
 \end{pmatrix}.
 \end{split}
\end{equation}
Its constrained inverse is $Q$, because
$A^{-1}\ee=\ee$ and $\ee\T A^{-1}\ee=5$.
The constrained inverse of $A_{II}$ for $I=\{1,2,3,4\}$ is $S$.

\subsection{Exact principal-inverse certificate}
Table~\ref{tab:certificate} lists
\[
 \mu(T)=\min_{i\in T}(H_T\ee_T)_i,
 \qquad \eta(T)=\min_{i\in T}\eta_i(T),\qquad H_T=Q_{TT}^{-1},
\]
for all $30$ nonempty proper subsets. Concatenated digits denote index
sets. All entries are rational and can be checked directly. In particular,
with $G=1058Q$ as in \eqref{eq:Qexplicit}, each inverse is given by
\[
 H_T=1058\,\frac{\operatorname{adj}(G_{TT})}{\det(G_{TT})}.
\]
Its row sums and negative-part row margins give the two entries for $T$.

\begin{table}[htbp]
\centering\small
\begin{tabular}{@{}rrr@{\qquad}rrr@{}}
\toprule
$T$&$\mu(T)$&$\eta(T)$&$T$&$\mu(T)$&$\eta(T)$\\
\midrule
1&$1058/1553$&$1058/1553$&123&$71/131$&$243/524$\\
2&$529/969$&$529/969$&124&$43/1432$&$43/1432$\\
3&$1058/3937$&$1058/3937$&125&$2/3$&$5/12$\\
4&$2/3$&$2/3$&134&$354/193$&$2289/3088$\\
5&$2$&$2$&135&$8$&$2$\\
12&$893/2433$&$893/2433$&145&$0$&$0$\\
13&$419/340$&$1553/2720$&234&$5731/1248$&$1225/1248$\\
14&$1760/2289$&$1553/2289$&235&$31/22$&$39/88$\\
15&$1$&$23/32$&245&$351/368$&$3/4$\\
23&$154/331$&$102/331$&345&$0$&$0$\\
24&$2921/1225$&$1587/1225$&1234&$1363/92$&$78/23$\\
25&$35/39$&$23/39$&1235&$5$&$1$\\
34&$46/119$&$1587/5593$&1245&$1/4$&$-1/2$\\
35&$2/3$&$23/72$&1345&$0$&$0$\\
45&$0$&$0$&2345&$197/92$&$1$\\
\bottomrule
\end{tabular}
\caption{Complete certificate for the matrix \eqref{eq:Qexplicit}.}
\label{tab:certificate}
\end{table}

Since every $\mu(T)$ is nonnegative, $Q$ has FMP by
Theorem~\ref{th:inverse}. There is exactly one negative row margin:
$T_0=\{1,2,4,5\}$, at row $1$. The complete exceptional inverse is
\begin{equation}\label{eq:exception}
 H_{T_0}=\begin{pmatrix}
 1&-3/4&-3/4&3/4\\
 -3/4&2&2&-5/4\\
 -3/4&2&3&-9/4\\
 3/4&-5/4&-9/4&179/46
 \end{pmatrix}.
\end{equation}
Its row sums are $(1/4,2,2,105/92)$ and its margins are
$(-1/2,0,0,9/23)$.

\begin{lemma}\label{lem:fourfacets}
The four facets obtained by deleting $h=1,2,3,4$ have hereditary
correctness. Every proper face of $I=\{1,2,3,4\}$ also has that property.
\end{lemma}
\begin{proof}
For the facet obtained by deleting $h$, formula \eqref{eq:subinverse}
says that the inverses used in its hereditary test are
$(H_T)_{T\setminus\{h\},T\setminus\{h\}}$ with $h\in T\subsetneq P$.
Deleting an off-diagonal entry cannot decrease a row margin of the form
\eqref{eq:eta}. Every original margin is nonnegative except the indicated
row of \eqref{eq:exception}.
Deleting $h=1$ removes that row. Deleting $h=2$ or $h=4$ removes one of
its entries $-3/4$, raising its margin to $1/4$.
The exceptional set does not contain $h=3$, so it does not occur in the
test for that facet. Thus all four facets satisfy
Theorem~\ref{th:hereditary}. Each proper face of $I$ omits some index
in $\{1,2,3,4\}$ and lies in one of these hereditary facets.
\end{proof}

\subsection{Failure on the remaining facet}
For $I=\{1,2,3,4\}$, take
\begin{equation}\label{eq:badb}
 \begin{split}
 x^*&=e_1,\qquad r=(0,9,0,12)\T,\\
 b_I&=A_{II}e_1-r
 =\left(\frac{1513}{1150},-\frac{5626}{575},
 \frac{841}{2300},-\frac{7351}{575}\right)\T.
 \end{split}
\end{equation}
These are KKT data with $\nu=0$, so $e_1$ is the true optimum.
The affine minimizer is
\begin{equation}\label{eq:bady}
 z=e_1-Sr=\frac1{23}(-85,-75,252,-69)\T.
\end{equation}
Coordinate $1$ is the unique negative minimum. Thus $\MIN$ immediately
removes a positive optimal coordinate, and hereditary correctness fails.

\subsection{Correctness on the reachable part of that facet}
For the full matrix \eqref{eq:Qexplicit}, direct inspection gives
\begin{equation}\label{eq:csets}
 C_1=C_2=C_3=C_4=\varnothing,\qquad C_5=\{4\}.
\end{equation}
Only facet $I$ requires further analysis. The condition inherited from
a first deletion of $5$ is $z_4\ge0$; the failing vector
\eqref{eq:bady} has $z_4=-3$.

\begin{lemma}\label{lem:guardeds}
On the facet with constrained inverse $S$, $\MIN$ is correct for every
linear term whose initial affine minimizer satisfies $z_4\ge0$, and for
every resolution of ties.
\end{lemma}
\begin{proof}
Let $x,r$ be the complementary KKT pair of an arbitrary such task;
then $z=x-Sr$. It suffices to prove that every first deletion is safe,
because Lemma~\ref{lem:fourfacets} covers all subsequent proper faces.

Rows $3$ and $4$ of $S$ have no positive off-diagonal entries. If $x_i>0$
for either index, then $r_i=0$ and $z_i\ge x_i>0$, so that index cannot
be removed.

If $x_1>0$ and $z_1<0$, then $r_1=0$, and the expression for $z_1$
forces $r_2>0$ and hence $x_2=0$. Therefore
\[
 z_1-z_2=x_1+\frac{27}{23}r_2+\frac{16}{23}r_3-r_4.
\]
If $r_4=0$, this is positive. Otherwise $x_4=0$ and the condition
$z_4=r_2-r_4\ge0$ gives $r_4\le r_2$. The difference is then at least
$x_1+4r_2/23+16r_3/23>0$. Thus coordinate $1$ cannot minimize $z$.

If $x_2>0$ and $z_2<0$, then $r_2=0$ and
\[
 12r_1>23x_2+28r_3+23r_4.
\]
In particular $r_1>0$ and $x_1=0$. Consequently
\[
 \begin{split}
 z_2-z_1
 &=x_2+\frac{20}{23}r_1-\frac{16}{23}r_3+r_4\\
 &>\frac83(x_2+r_4)+\frac43r_3>0.
 \end{split}
\]
Coordinate $2$ cannot minimize $z$ either. These strict comparisons
exclude every positive optimal coordinate from every possible first
choice, including all ties between the remaining indices.
\end{proof}

\begin{proposition}\label{prop:counter}
The matrix $A$ in \eqref{eq:Aexplicit} has full-start correctness and
does not have hereditary correctness.
\end{proposition}
\begin{proof}
Table~\ref{tab:certificate} proves FMP for $Q$. The four facets with
$h\le4$ are hereditary by Lemma~\ref{lem:fourfacets}. For $h=5$, the
condition prescribed by \eqref{eq:csets} is exactly $z_4\ge0$, under
which Lemma~\ref{lem:guardeds} applies. Theorem~\ref{th:guarded}
therefore gives full-start correctness. Equations
\eqref{eq:badb}--\eqref{eq:bady} prove failure on a face.
\end{proof}

\section{Minimality and the generic contrast}\label{sec:small}
On a two-coordinate simplex, every centered positive semidefinite
constrained inverse has the form
$a\left(\begin{smallmatrix}1&-1\\-1&1\end{smallmatrix}\right)$ with $a>0$;
thus its first step is safe. A singleton is trivial.

\begin{lemma}\label{lem:triangle}
For three coordinates, full-start correctness is equivalent to
\begin{equation}\label{eq:triangle}
 Q_{ij}\le\min(Q_{ii},Q_{jj})\qquad(i\ne j).
\end{equation}
\end{lemma}
\begin{proof}
For a positive definite $2\times2$ block
$\left(\begin{smallmatrix}a&c\\c&d\end{smallmatrix}\right)$,
the inverse applied to $(1,1)\T$ is $(d-c,a-c)\T/(ad-c^2)$.
Thus \eqref{eq:triangle} is exactly Theorem~\ref{th:inverse} for a
three-coordinate constrained inverse. It characterizes safety of the
first step. Subsequent faces have at most two coordinates and are
automatically correct.
\end{proof}

\begin{theorem}\label{th:four}
For $n\le4$, full-start correctness implies hereditary correctness,
without strictness assumptions and with all ties allowed.
\end{theorem}
\begin{proof}
The cases $n\le3$ follow from the preceding observations. Let $n=4$ and
suppose full-start correctness holds. Then $Q$ has FMP. Suppose a
three-coordinate facet $I=P\setminus\{h\}$ is not universally correct,
and write $S=Q_I$. By Lemma~\ref{lem:triangle}, there are distinct
$i,j\in I$ with $S_{ij}>S_{jj}>0$. Denote the third index by $k$.
Since $S\ee=0$, we have $S_{kj}=-S_{ij}-S_{jj}<0$.

On $I$, choose optimum $e_i$ and residual $te_j$ with
$t>1/(S_{ij}-S_{jj})$. These data determine a linear term on the facet,
and its affine minimizer $z=e_i-tSe_j$ satisfies
\[
 z_i<z_j<0<z_k.
\]
Put $q=Q_{hh}$ and $c_\ell=Q_{\ell h}/q\le1$.
We claim $c_j<1$. If $c_j=1$, the column inequality $Q_{ij}\le Q_{jj}$,
together with
$Q_{ij}=S_{ij}+qc_i$ and $Q_{jj}=S_{jj}+q$,
gives $c_i<1$. In the proper block $R=\{i,j,h\}$, block inversion
expresses the first two components of $Q_{RR}^{-1}\ee_R$ as
\[
 S_{\{i,j\},\{i,j\}}^{-1}(1-c_i,1-c_j)\T.
\]
Its $j$ component is
$-S_{ij}(1-c_i)/(S_{ii}S_{jj}-S_{ij}^2)<0$,
contradicting Theorem~\ref{th:inverse}. This proves the claim.

Use reverse insertion with $y_h=-L$ and
$y_\ell=z_\ell-Lc_\ell$, for large $L>0$. The optimum remains $e_i$.
If $c_i<1$, both $i$ and $j$ lie strictly above $h$ for sufficiently
large $L$. So does $k$: if $c_k<1$ this follows by growth of the
difference, while if $c_k=1$ the difference is $z_k>0$.
Thus $h$ is uniquely deleted first, and $i$ is then wrongly deleted on
the facet.
If $c_i=1$, then $y_i-y_h=z_i<0$. For large $L$, coordinate $i$ lies
strictly below $j$ because $c_j<1$, and strictly below $k$ because either
$c_k<1$ or $z_k>z_i$. In this case $i$ is wrongly deleted immediately.
Both alternatives contradict full-start correctness.
Hence all three-coordinate facets are universally correct, and all
smaller faces are automatically correct.
\end{proof}

Theorem~\ref{th:main} now follows from Proposition~\ref{prop:counter}
and Theorem~\ref{th:four}. In geometric terminology, the smallest
simplex dimension for this distinction is four, corresponding to five
vertices.

\begin{theorem}\label{th:generic}
Suppose that for every face $J$ with $|J|\ge2$ and distinct $i,j\in J$,
\begin{equation}\label{eq:generic}
 (Q_J)_{ij}\ne(Q_J)_{jj}.
\end{equation}
Then full-start correctness and hereditary correctness are equivalent.
The matrices satisfying \eqref{eq:generic} form an open dense subset
of the positive definite cone.
\end{theorem}
\begin{proof}
Universal correctness on a given face requires FMP, hence diagonal
maxima in each column. Under \eqref{eq:generic}, those maxima are strict.
Every equality set in Lemma~\ref{lem:reach} is therefore empty.
Any failing child task could be extended by a first deletion, contradicting
correctness on the parent. Induction gives correctness on all faces.
The reverse implication is immediate.

The excluded equalities are finitely many rational equations in the
entries of $A$, with denominators nonzero on the positive definite cone.
None is identically zero: for $A=I$, every face has
$Q_J=I_J-\ee_J\ee_J\T/|J|$, with a diagonal entry one larger than any
off-diagonal entry. Clearing denominators yields finitely many proper
algebraic zero sets, each with empty interior. Their complement is open
and dense.
\end{proof}

\begin{corollary}
For $n=5$, full-start correctness implies hereditary correctness whenever
each diagonal entry of the full $Q$ strictly exceeds the other entries
in its column. No generic assumption on smaller faces is needed.
\end{corollary}
\begin{proof}
Reverse insertion gives universal correctness on every four-coordinate
facet, and Theorem~\ref{th:four} makes each such facet hereditary.
\end{proof}

The equality $Q_{45}=Q_{55}$ in \eqref{eq:Qexplicit} is thus essential
to the exhibited distinction. It imposes a restriction on the tasks
that can follow a first deletion, although the conclusion of full-start
correctness itself quantifies over every linear term and every tie.

\section{Recognition of the first-step condition and a limit of shorter tests}
\label{sec:recognition}

The principal-submatrix descriptions above do not automatically require
exponential recognition time. In fact, the first-step maximum principle
has a polynomial-time test, including its boundary. This does not give
such a test for hereditary correctness.
Throughout this section, complexity is measured in the total binary
encoding length of the rational input, using the standard Turing model.

\begin{lemma}[Positive diagonal regularization]\label{lem:regularize}
Let $M$ be SPD and $d>0$. Suppose $M_{TT}^{-1}d_T\ge0$ for every
nonempty principal set $T$. If $D$ is any positive diagonal matrix, then
$(M+D)_{TT}^{-1}d_T>0$ for every such $T$.
\end{lemma}
\begin{proof}
Write $v^S=M_{SS}^{-1}d_S$. Multilinearity in the diagonal entries of
Cramer's numerator gives, for $i\in T$,
\[
 \det(M_{TT}+D_{TT})\bigl[(M_{TT}+D_{TT})^{-1}d_T\bigr]_i
 =\sum_{i\in S\subseteq T}\det(M_{SS})v_i^S
       \prod_{j\in T\setminus S}D_{jj}.
\]
All summands are nonnegative. The term $S=\{i\}$ equals
$d_i\prod_{j\in T\setminus\{i\}}D_{jj}>0$. The determinant on the
left is positive, proving the assertion.
\end{proof}

\begin{theorem}[Polynomial recognition of FMP]\label{th:fmp-poly}
For rational SPD $M$, the condition $M_{TT}^{-1}\ee_T\ge0$ on all
nonempty principal sets is recognizable in polynomial bit time.
For rational symmetric $Q\succeq0$ with
$\ker Q=\operatorname{span}(\ee)$, FMP is likewise recognizable in
polynomial bit time. In particular, safety of the first $\MIN$ step for
every linear term has a polynomial-time recognition algorithm.
\end{theorem}
\begin{proof}
We use the following known strict criterion. For a P-matrix $N$ and
$d>0$, Morris's Theorem~6.1 \cite{morris} states that
$N_{TT}^{-1}d_T>0$ on every nonempty principal set if and only if there
is a Z-matrix $X$ such that $N\T X$ is also a Z-matrix and
$X\T d>0$. Here a Z-matrix has nonpositive off-diagonal entries.
For fixed rational $N$ and $d=\ee$, this is the linear feasibility
problem
\begin{equation}\label{eq:morris-lp}
 X_{ij}\le0,\quad (N\T X)_{ij}\le0\quad(i\ne j),
 \qquad X\T\ee\ge\ee.
\end{equation}
Homogeneous scaling justifies the last normalization. There are $n^2$
variables and polynomially many inequalities. An SPD matrix satisfies
the P-matrix hypothesis.

To reduce the weak condition to this strict test, let $L$ be the total
binary encoding length of the input matrix and set, for example,
\[
 h=4n\bigl(L+\lceil\log_2(n+1)\rceil+1\bigr),
 \qquad \varepsilon=2^{-5h-2}.
\]
Uniformly over the SPD principal blocks, determinant and eigenvalue
bounds give
\begin{equation}\label{eq:rational-gap}
 \lambda_{\min}(M_{TT})\ge2^{-h},\qquad
 \bigl|(M_{TT}^{-1}\ee_T)_i\bigr|\ge2^{-h}
 \quad\hbox{whenever that component is nonzero}.
\end{equation}
For clarity, a common denominator is at most $2^L$, each entry has
absolute value at most $2^L$, every nonzero determinant of order at
most $n$ has magnitude at least $2^{-nL}$, and Hadamard's inequality
bounds such determinants above by $(n2^L)^n$. These estimates, Cramer's
rule, and $\lambda_{\max}\le n2^L$ imply
\eqref{eq:rational-gap} with the stated conservative $h$.

The resolvent identity now yields, on every principal block,
\[
 \left\|(M_{TT}+\varepsilon I)^{-1}\ee_T-M_{TT}^{-1}\ee_T\right\|_2
 \le\varepsilon\,2^{2h}\sqrt n<2^{-h-1}.
\]
Thus a negative component remains negative. If all components are
nonnegative, Lemma~\ref{lem:regularize} makes all of them strictly
positive. Consequently weak recognition for $M$ is exactly strict
recognition for $N=M+\varepsilon I$, which is decided by
\eqref{eq:morris-lp}. The bit length of $\varepsilon$ is polynomial.

For centered $Q$, apply the same bounds to its proper principal blocks
and use $N=Q+\varepsilon I$. The full-block solution is automatically
$N^{-1}\ee=\varepsilon^{-1}\ee>0$. For every proper block the preceding
argument applies unchanged. Theorem~\ref{th:inverse} then identifies
the tested condition with FMP. The one-coordinate case is automatic.
For rational SPD input $A$, the constrained inverse $Q$ is computable by
exact rational linear algebra with polynomial encoding length, so the
same conclusion applies to first-step recognition from $A$.
\end{proof}

The LP ingredient is a classical strict $n$-step-vector recognition
theorem. The diagonal argument above supplies the non-strict extension
used here. Weak principal-inverse inequalities are themselves classical:
Chu~\cite[p.~463, equation (2)]{chu} attributes them to
Pang--Chandrasekaran. Chu also notes there that, even for P-matrices,
the weak-vector set need not be the closure of the strict-vector set.
Our proof does not use that false closure inference: it perturbs the
matrix and controls the negative components by a rational separation
bound. Chu's algorithms take such a vector as given;
the discussion on p.~469 leaves efficient checking of the broader
extended conditions as a separate question. Theorem~\ref{th:fmp-poly}
is an explicit deduction for the SPD and centered cases; it does not
resolve recognition for all extended $n$-step conditions. In particular,
a centered $Q$ is not being asserted to have an extended $n$-step vector
in Chu's sense: its full singular block cannot solve $Qx=\ee$.
The test here involves only its proper principal blocks. No new
definition of an $n$-step vector or complementary pivot method is needed.
Testing FMP only for the full $Q$ does not test it for all
Schur complements $Q_J$; the counterexample in
Section~\ref{sec:example} already separates those conditions.

One might try to supplement the full FMP test with the row conditions
from Theorem~\ref{th:hereditary} on only the $n$ maximal proper blocks.
The next example rules out even that stronger shortcut.

\begin{proposition}\label{prop:recognition-shortcut}
Full FMP together with
$\eta_i(P\setminus\{k\})\ge0$ for all $k$ and all $i\ne k$ is
equivalent to hereditary correctness when $n\le4$.
For $n=5$ these tests can all hold strictly while hereditary correctness
fails. This failure persists on an open subset of the SPD cone.
\end{proposition}
\begin{proof}
For $n\le4$, every proper block not covered by a maximal-block test
has size at most two. A one-dimensional SPD inverse has positive
margin. In a two-dimensional SPD inverse, a negative row margin
requires a negative off-diagonal entry; that margin is then the full
row sum, already excluded by FMP and Theorem~\ref{th:inverse}.
Theorem~\ref{th:hereditary} proves the equivalence in these dimensions.

For the five-coordinate construction, put
\begin{equation}\label{eq:recognition-construction}
 \begin{split}
 D&=\begin{pmatrix}
 1&-6/5&1/2\\[-2pt]
 -6/5&23/10&-1\\[-2pt]
 1/2&-1&1
 \end{pmatrix},\qquad w=\tfrac12\ee_3,\\
 H&=\begin{pmatrix}D+ww\T&w\\w\T&1\end{pmatrix},
 \qquad
 A=\begin{pmatrix}H+\ee_4\ee_4\T&\ee_4\\\ee_4\T&1\end{pmatrix}.
 \end{split}
\end{equation}
The leading principal minors of $D$ are $1,43/50,97/200$.
Two Schur-complement identities show that $H$ and $A$ are SPD.
Let $Q$ be the constrained inverse of $A$. Since $H$ is the edge Gram
matrix anchored at coordinate $5$,
$Q_{1234,1234}^{-1}=H$ and
\[
 Q_{123,123}^{-1}=H/H_{44}=D,
 \qquad (\eta_1(123),\eta_2(123),\eta_3(123))=(-1/5,1/10,0).
\]
Thus hereditary correctness fails.

Table~\ref{tab:recognition-shortcut} gives all proper-block row-sum
and margin minima. Every row-sum minimum is strictly positive, so
$Q$ has FMP. Every maximal proper block has strictly positive margin
minimum. For completeness, their full margin vectors, in the order of
the remaining indices, are
\[
\begin{array}{c|c}
\text{omitted index}&\text{margin vector}\\\hline
1&(57/10,1,5/4,5/4)\\
2&(57/10,53/10,51/20,51/20)\\
3&(1,53/10,5/4,5/4)\\
4&(3/10,17/20,1/2,1)\\
5&(3/10,17/20,1/2,1).
\end{array}
\]
All asserted good tests are strict, and the failing margin is strictly
negative. There are finitely many tests, all continuous in $A$ on the
SPD cone. The same separation therefore holds on an open neighborhood
of this matrix.
\end{proof}

\begin{table}[htbp]
\centering\small
\begin{tabular}{@{}rrr@{\qquad}rrr@{}}
\toprule
$T$&$\mu(T)$&$\eta(T)$&$T$&$\mu(T)$&$\eta(T)$\\
\midrule
1&$97/260$&$97/260$&123&$1/10$&$-1/5$\\
2&$97/150$&$97/150$&124&$1/2$&$3/10$\\
3&$97/172$&$97/172$&125&$1/2$&$3/10$\\
4&$194/715$&$194/715$&134&$349/170$&$457/510$\\
5&$194/715$&$194/715$&135&$349/170$&$457/510$\\
12&$1/20$&$1/20$&145&$93/106$&$26/53$\\
13&$81/230$&$81/230$&234&$41/50$&$31/50$\\
14&$31/35$&$52/105$&235&$41/50$&$31/50$\\
15&$31/35$&$52/105$&245&$67/20$&$1$\\
23&$7/20$&$7/20$&345&$101/190$&$86/285$\\
24&$67/40$&$3/4$&1234&$27/20$&$3/10$\\
25&$67/40$&$3/4$&1235&$27/20$&$3/10$\\
34&$303/457$&$172/457$&1245&$23/10$&$1$\\
35&$303/457$&$172/457$&1345&$57/5$&$51/20$\\
45&$97/521$&$97/521$&2345&$27/10$&$1$\\
\bottomrule
\end{tabular}
\caption{Complete rational certificate for
\eqref{eq:recognition-construction}. As before,
$\mu(T)=\min_i(Q_{TT}^{-1}\ee_T)_i$ and
$\eta(T)=\min_i\eta_i(T)$.}
\label{tab:recognition-shortcut}
\end{table}

Hereditary recognition belongs to coNP: a proper set $T$ and a row with
$\eta_i(T)<0$ are a polynomial-size certificate. The polynomial FMP
algorithm and the failed shortcut alone do not classify the stronger
recognition problem. Section~\ref{sec:hereditary-hardness} establishes its
coNP-completeness by a separate reduction.

\section{Hereditary recognition is coNP-complete}
\label{sec:hereditary-hardness}

The polynomial first-step test does not extend to hereditary correctness
unless $\mathrm P=\mathrm{NP}$. The difficulty persists even when every
inequality in a positive instance of the recognition problem is strict.

\begin{theorem}\label{th:hereditary-hardness}
Given a rational symmetric positive definite matrix $A$, recognizing
hereditary correctness of $\MIN$ is coNP-complete under polynomial-time
many-one reductions. The reduction has strictly positive principal-inverse
row margins on its hereditary instances and a strictly negative margin on
its other instances.
\end{theorem}

We prove the theorem using a Gram construction. The construction must
control every proper inverse block of the centered matrix $Q$, including
blocks containing the coordinate chosen as the initial affine anchor.
One selected row margin encodes the graph's clique number. The estimates
then show that every other row and every change of anchor remain valid
when the graph has no clique of the requested size. Keeping these two
checks explicit is necessary: a reduction for a single anchored block
would not prove hereditary hardness.
The following identity separates that issue from ordinary Schur complements.
For any symmetric matrix $H$, write
\[
 \eta_i(H)=H_{ii}+\sum_{j\ne i}\min(H_{ij},0).
\]
We use $G/G_{EE}$ for the Schur complement obtained by eliminating $E$;
when $E=\varnothing$, this means $G$ itself. Empty products and sums
below have their usual zero-dimensional meanings.

\begin{lemma}[Changing the omitted anchor]\label{lem:anchor-lift}
Let $G\succ0$ have index set $V$, introduce a coordinate $0$, and set
\[
 A=\begin{pmatrix}G+\ee\ee\T&\ee\\\ee\T&1\end{pmatrix},
 \qquad Q=\begin{pmatrix}I\\-\ee\T\end{pmatrix}
              G^{-1}\begin{pmatrix}I&-\ee\end{pmatrix}.
\]
Then $A$ is SPD and $Q$ is its constrained inverse. If $0\notin T$, then
$Q_{TT}^{-1}=G/G_{EE}$, where $E=V\setminus T$.
If $T=R\cup\{0\}$ is proper, put $O=V\setminus R\ne\varnothing$ and
\[
 S=G/G_{OO},\qquad h=\ee_R-G_{RO}G_{OO}^{-1}\ee_O,
 \qquad \beta=(\ee_O\T G_{OO}^{-1}\ee_O)^{-1}.
\]
In the order $R,0$, one has
\begin{equation}\label{eq:anchor-lift}
 Q_{TT}^{-1}=
 \begin{pmatrix}S+\beta hh\T&\beta h\\\beta h\T&\beta\end{pmatrix}.
\end{equation}
Consequently, if every principal Schur complement of $G$ has nonnegative
row margins and
\begin{equation}\label{eq:anchor-h-positive}
 \ee_R-G_{RO}G_{OO}^{-1}\ee_O\ge0
 \quad\text{for every nonempty }O\subseteq V,
 \quad R=V\setminus O,
\end{equation}
then $A$ has hereditary correctness.
\end{lemma}
\begin{proof}
The Schur complement of the last entry of $A$ is $G$. Direct substitution
in \eqref{eq:qdef} gives $Q$. The first inverse identity is the usual
principal-inverse Schur identity. For the second, block inversion gives
$Q_{RR}=S^{-1}$, $Q_{R0}=-S^{-1}h$, and
$Q_{00}=h\T S^{-1}h+\beta^{-1}$, which yields
\eqref{eq:anchor-lift}. Its added rank-one matrix is entrywise nonnegative
when $h\ge0$. Increasing a diagonal or an off-diagonal entry cannot decrease
$\eta_i$. Thus the margins in the rows of $R$ are at least those in $S$,
and the anchor row has margin $\beta>0$. Theorem~\ref{th:hereditary}
applies. If $R$ is empty, \eqref{eq:anchor-lift} is the scalar $\beta$.
\end{proof}

\subsection{The reduction and its target margin}
Take a CLIQUE instance \cite{karp} consisting of a simple graph $\mathcal G$ on $m\ge2$
vertices and an integer $2\le k\le m$. Define the symmetric matrix $W$ by
\[
 W_{ii}=1,\qquad
 W_{ij}=\begin{cases}0,&\{i,j\}\in E(\mathcal G),\\
 -1,&\{i,j\}\notin E(\mathcal G),\ i\ne j.
 \end{cases}
\]
This signed matrix also occurs in the biclique reduction of
Tan~\cite[Lemma~4]{tan}; the weighted score used below has a different role.
The auxiliary set $C$ consists of left and right copies of these vertices.
Put
\begin{equation}\label{eq:hard-parameters}
 \begin{gathered}
 J=\begin{pmatrix}0&W\\W&0\end{pmatrix},\quad
 D=I_{2m}+\varepsilon J,\quad
 a=\binom{\ee_m}{0},\quad b=\binom{0}{\ee_m},\quad d=a+b,\\
 r=8m,\quad t=\frac1{100m},\quad s=\frac tr=\frac1{800m^2},\quad
 \varepsilon=\frac1{128m^3},\quad
 \delta=r\varepsilon(k-\tfrac12).
 \end{gathered}
\end{equation}
In particular $\varepsilon m\le1/512$ and $0<\delta<1/(16m)$.
On $V=\{u,v\}\cup C$, define the Gram matrix
\begin{equation}\label{eq:hard-gram}
 G=\begin{pmatrix}
 1&-1+s^2\delta&st a\T\\
 -1+s^2\delta&2-2s^2\delta&-st d\T\\
 st a&-st d&t^2D
 \end{pmatrix}.
\end{equation}
The output matrix is the $A$ in Lemma~\ref{lem:anchor-lift}.

To verify positive definiteness, consider
\begin{equation}\label{eq:hard-original-gram}
 M=\begin{pmatrix}
 1&-s^2\delta&s a\T\\
 -s^2\delta&1&s b\T\\
 s a&s b&D
 \end{pmatrix}.
\end{equation}
An endpoint row has off-diagonal absolute sum $ms+s^2\delta<1/4$;
an auxiliary row has such sum at most $s+\varepsilon m<1/4$.
Hence $M\succeq(3/4)I$. If its Gram vectors are $p,q,c$, the invertible
change $u=p$, $v=-p-q$, $c_{\rm new}=tc$ gives \eqref{eq:hard-gram}.
Thus $G$ and the output $A$ are SPD.

For left and right subsets $L,R\subseteq[m]$, define
\[
 q_j(L)=\sum_{i\in L}W_{ij},\quad
 F(L,R)=\sum_{j\in R}q_j(L),\quad
 P(L,R)=\sum_{j\notin R}\max(q_j(L),0).
\]
Each $q_j(L)$ is an integer at most one. It is positive precisely when
$j\in L$ and $j$ is adjacent to every other member of $L$, in which case
it equals one. These positive columns form a clique $B(L)$. Therefore
\begin{equation}\label{eq:hard-score}
 F(L,R)+rP(L,R)
 \le |B(L)\cap R|+r|B(L)\setminus R|
 \le r\,\omega(\mathcal G).
\end{equation}
For a $k$-clique $L$ and $R=\varnothing$, the score equals $rk$.

Eliminate an arbitrary auxiliary set $E\subseteq C$, whose left and right
parts are $L,R$. The residual Gram of $p,q$ in \eqref{eq:hard-original-gram}
is
\[
 F_E=\begin{pmatrix}A_E&K_E\\K_E&B_E\end{pmatrix},\qquad
 \begin{aligned}
 A_E&=1-s^2a_E\T D_{EE}^{-1}a_E,\\
 B_E&=1-s^2b_E\T D_{EE}^{-1}b_E,\\
 K_E&=-s^2(\delta+f_E),\quad f_E=a_E\T D_{EE}^{-1}b_E.
 \end{aligned}
\]
Empty sets have their usual zero-product interpretation. A Schur complement
preserves a uniform lower bound $M\succeq cI$: its quadratic form at $y$
is the minimum of $(x,y)\T M(x,y)$ over $x$, hence is at least
$c\|y\|^2$. Uniformly in $E$,
\begin{equation}\label{eq:hard-endpoint-bounds}
 F_E\succeq\tfrac34I,\quad A_E,B_E\le1,\quad
 |K_E|\le s^2(\delta+2m)<\tfrac1{16}.
\end{equation}
For $c\in C\setminus E$, put
\[
 \alpha_c=a_c-D_{cE}D_{EE}^{-1}a_E,\qquad
 \gamma_c=1-D_{cE}D_{EE}^{-1}\ee_E.
\]
The residual Gram of $u,v$ is
\begin{equation}\label{eq:hard-U}
 U=\begin{pmatrix}A_E&-A_E-K_E\\-A_E-K_E&A_E+2K_E+B_E\end{pmatrix};
\end{equation}
its coupling to $c$ is $st(\alpha_c,-\gamma_c)\T$, and the auxiliary
block is $t^2\bar D$, where $\bar D=D/D_{EE}$.
Since $U_{uv}<0$, the target row margin is exactly
\begin{equation}\label{eq:hard-target}
 \eta_u=s^2(\delta+f_E)+st\sum_{c\in C\setminus E}\min(\alpha_c,0).
\end{equation}

The Neumann expansion, $\|J\|_2\le m$, and
$\|D_{EE}^{-1}\|_2\le2$ give
\[
 |f_E+\varepsilon F(L,R)|\le2\varepsilon^2m^3,
 \qquad
 |\alpha_c-(a_c-\varepsilon a_E\T J_{Ec})|
 \le2\varepsilon^2m^2.
\]
For the second estimate, use $\|a_E\|_2\le\sqrt m$,
$\|D_{EE}^{-1}-I\|_2\le2\varepsilon m$, and
$\|D_{Ec}\|_2\le\varepsilon\sqrt m$.
The expression approximating $\alpha_c$ equals one on a remaining left
index and $-\varepsilon q_c(L)$ on a remaining right index.
The function $\min(x,0)$ is 1-Lipschitz, so \eqref{eq:hard-target} yields
\begin{equation}\label{eq:hard-target-error}
 \left|\frac{\eta_u}{s^2}
       -\bigl[\delta-\varepsilon(F(L,R)+rP(L,R))\bigr]\right|
 \le6r\varepsilon^2m^3<\frac{r\varepsilon}{16}.
\end{equation}
If $\omega(\mathcal G)\le k-1$, equations
\eqref{eq:hard-score}--\eqref{eq:hard-target-error} imply
\begin{equation}\label{eq:hard-no-margin}
 \eta_u>\frac7{16}s^2r\varepsilon>0
 \quad\text{for every }E\subseteq C.
\end{equation}
If $L$ is a $k$-clique, eliminate exactly its left copy. Then $D_{EE}=I$,
$f_E=0$, and precisely the $k$ corresponding right indices have negative
$\alpha_c=-\varepsilon$. Thus the exact margin is
\begin{equation}\label{eq:hard-yes-margin}
 \eta_u=s^2(\delta-r\varepsilon k)
       =-\tfrac12s^2r\varepsilon<0.
\end{equation}

\subsection{All other rows and all omitted anchors}
We must exclude unintended negative margins when the graph has no
$k$-clique. First, uniformly in the auxiliary set $E$,
\begin{equation}\label{eq:hard-aux-bounds}
 \begin{gathered}
 \|D_{EE}^{-1}\|_\infty\le2,\quad
 |\alpha_c|,|\gamma_c|\le2,\quad
 \bar D_{cc}\ge1-2\varepsilon^2m^2,\\
 \sum_{j\ne c}|\bar D_{cj}|\le2\varepsilon m.
 \end{gathered}
\end{equation}
Indeed the absolute row sums of $\varepsilon J$ are at most
$\varepsilon m$, and those of the Schur correction are at most
$2\varepsilon^2m^2$. When neither endpoint is eliminated, the $v$-row
margin is at least
\[
 B_E+K_E-st\sum_c|\gamma_c|\ge\frac{11}{16}-4mst>0.
\]
Each auxiliary row has margin at least
\begin{equation}\label{eq:hard-aux-margin}
 t^2\bigl(1-2\varepsilon^2m^2-2\varepsilon m-4/r\bigr)>0.
\end{equation}

Next eliminate a nonempty subset $Z\subseteq\{u,v\}$, after $E$.
The smallest squared singular value of the change from $p,q$ to $u,v$
is $(3-\sqrt5)/2>1/3$. Equations
\eqref{eq:hard-endpoint-bounds}--\eqref{eq:hard-U} give $U\succeq I/4$,
so $\|U_{ZZ}^{-1}\|_2\le4$. Each auxiliary-to-endpoint row has norm at
most $2\sqrt2st$. Every entry of this additional Schur correction has
modulus at most $32s^2t^2$. Consequently the new auxiliary diagonal is
at least $3t^2/4$, and its off-diagonal absolute row sum is at most
\[
 t^2(2\varepsilon m+64ms^2).
\]
If one endpoint remains, its diagonal is at least $1/4$. Its coupling
to each auxiliary has modulus at most $5st$, because the endpoint
elimination multiplier satisfies
\[
 \frac{|U_{uv}|}{U_{uu}}\le\frac{17}{12}<\frac32,
 \qquad
 \frac{|U_{uv}|}{U_{vv}}\le\frac{17/16}{11/8}<1.
\]
The remaining endpoint margin is therefore at least $1/4-10mst>0$.
Every auxiliary margin is at least
\begin{equation}\label{eq:hard-postendpoint-margin}
 t^2\bigl(3/4-5/r-2\varepsilon m-64ms^2\bigr)>0.
\end{equation}
All displayed positive bounds follow directly from $m\ge2$ and
\eqref{eq:hard-parameters}; for example $5/r\le5/16$,
$2\varepsilon m\le1/256$, and $64ms^2\le1/80000$.
Schur associativity allows every eliminated set to be treated in this
order. Together with \eqref{eq:hard-no-margin}, this controls all
principal Schur complements of $G$ in a no-clique instance.

It remains to prove \eqref{eq:anchor-h-positive}. This property in fact
holds for either CLIQUE answer. Fix a nonempty omitted set $O\subseteq V$,
and write $B=O\cap C$, $Z=O\cap\{u,v\}$. After eliminating $B$, the
endpoint right-hand side for $h$ is
\[
 f_u=1-r^{-1}a_B\T D_{BB}^{-1}\ee_B,\qquad
 f_v=1+r^{-1}\ee_B\T D_{BB}^{-1}\ee_B.
\]
Equation~\eqref{eq:hard-aux-bounds} implies
\begin{equation}\label{eq:hard-h-rhs}
 f_u\ge3/4,\quad f_v\ge1,\quad \|f\|_\infty\le3/2.
\end{equation}
Let $K$ be the endpoint residual after eliminating $B$. It has negative
off-diagonal entry by \eqref{eq:hard-endpoint-bounds}, and
\[
 K=K_0+E_0,\qquad
 K_0=\begin{pmatrix}1&-1\\-1&2\end{pmatrix},\qquad
 \|E_0\|_\infty\le8ms^2.
\]
The inverse of every nonempty principal block $K_{ZZ}$ has infinity norm
at most four. For the full block, this follows from
\[
 \|K^{-1}\|_\infty\le
 \frac{\|K_0^{-1}\|_\infty}{1-\|K_0^{-1}\|_\infty\|E_0\|_\infty}
 \le\frac3{1-24ms^2}<4;
\]
for a singleton it follows from $K_{uu}\ge3/4$ and $K_{vv}\ge11/8$.

For a remaining auxiliary $j$, put
$b_j=1-D_{jB}D_{BB}^{-1}\ee_B$. Then $3/4\le b_j\le5/4$, and
its endpoint coupling is $st(\alpha_j,-b_j)$ with $|\alpha_j|\le5/4$.
After eliminating $Z$, its component in \eqref{eq:anchor-h-positive} is
\[
 h_j=b_j-st(\alpha_j,-b_j)_ZK_{ZZ}^{-1}f_Z
 \ge3/4-15st>1/2.
\]
If an endpoint remains and $Z$ is empty, its component is its positive
$f$ component. Otherwise its component is
$f_i-K_{iZ}K_{ZZ}^{-1}f_Z>0$, since $Z$ is a singleton,
$K_{iZ}<0$, and $f_Z>0$. This proves
\eqref{eq:anchor-h-positive} for every omitted set, including an empty
$B$ or $Z$.

\begin{proof}[Completion of the proof of Theorem~\ref{th:hereditary-hardness}]
The construction has $2m+3$ coordinates and polynomially encoded rational
entries. If $\mathcal G$ has no $k$-clique, all principal Schur margins
of $G$ are strictly positive and \eqref{eq:anchor-h-positive} holds.
Lemma~\ref{lem:anchor-lift} makes every proper inverse margin of $Q$
strictly positive. If $\mathcal G$ has a $k$-clique,
\eqref{eq:hard-yes-margin} is a negative margin of a proper inverse
block omitting coordinate $0$. Theorem~\ref{th:hereditary} therefore
reduces CLIQUE to failure of hereditary correctness.
Finally, a proper set $T$ and a row $i$ with $\eta_i(T)<0$ form a
polynomial-size certificate, verifiable by exact rational linear algebra.
This proves coNP membership and completes the reduction.
\end{proof}

The large ratio $t/s$ has two roles: it protects the auxiliary rows and
the other affine anchors, while replacing the bilinear score $F$ by
$F+rP$, whose maximum remains exactly $r\omega(\mathcal G)$.
The same family also determines the complexity of full-start recognition.

\subsection{Full-start recognition}
\begin{theorem}\label{th:fullstart-hardness}
Given a rational SPD matrix $A$, recognizing full-start correctness of
$\MIN$ is coNP-complete. The hardness holds on matrices satisfying the
strict diagonal-maxima inequalities
\[
 (Q_J)_{ij}<\min\bigl((Q_J)_{ii},(Q_J)_{jj}\bigr)
 \quad\text{for every face }J\text{ and distinct }i,j\in J.
\]
\end{theorem}
\begin{proof}
We first show that every proper inverse block $H_T=Q_{TT}^{-1}$ of every
constructed instance, with either CLIQUE answer, satisfies
\begin{equation}\label{eq:hard-pair-margin}
 (H_T)_{ii}+(H_T)_{ij}>0\qquad(i\ne j).
\end{equation}
For an ordinary Schur complement of $G$, all rows except possibly $u$
when both endpoints remain have positive $\eta_i$, by the preceding
bounds independently of the graph answer. Such rows satisfy
\eqref{eq:hard-pair-margin}. For the exceptional row, the sum with the
$v$ entry is $s^2(\delta+f_E)$. Since $F(L,R)\le m$,
\[
 \delta+f_E\ge12m\varepsilon-m\varepsilon-2\varepsilon^2m^3>0.
\]
Its sum with an auxiliary entry is at least $3/4-2st>0$.
For blocks containing coordinate $0$, \eqref{eq:anchor-lift} preserves
the pair sums between ordinary indices by adding an entrywise
nonnegative matrix. The sums involving the anchor are
$S_{ii}+\beta h_i^2+\beta h_i>0$ and $\beta(1+h_i)>0$.
Thus \eqref{eq:hard-pair-margin} holds for every proper inverse block.

For a face $J$ with $|J|\ge3$, set $D=P\setminus J$ and
$T=D\cup\{i,j\}\subsetneq P$. Equation~\eqref{eq:subinverse} gives
\[
 \bigl((Q_J)_{\{i,j\},\{i,j\}}\bigr)^{-1}
       =(H_T)_{\{i,j\},\{i,j\}}.
\]
The two row sums on the right are positive by
\eqref{eq:hard-pair-margin}. Inverting a $2\times2$ SPD matrix shows
that these inequalities are exactly the asserted strict diagonal maxima
of $Q_J$. Two-coordinate faces satisfy them automatically.
Theorem~\ref{th:generic} therefore makes full-start and hereditary
correctness equivalent on every instance of the reduction. This proves
coNP-hardness of full-start recognition.

For membership, a failed run must delete a positive coordinate of its
unique simplex optimum. Otherwise that optimum remains feasible on the
terminal face, whose feasible affine minimizer must equal it. Fix a
permitted deletion prefix ending with such a deletion $h$, and fix an
optimal support $S$. The existence of a witness is a rational linear
system in $b,x,r,\nu$ and the affine minimizers $y^\ell$ with multipliers
$\lambda_\ell$ on the prefix faces $J_\ell$:
\[
 \begin{gathered}
 Ax-b-\nu\ee=r,\quad x,r\ge0,\quad\ee\T x=1,
 \quad x_{S^c}=0,\quad r_S=0,\\
 A_{J_\ell J_\ell}y^\ell-b_{J_\ell}
       =\lambda_\ell\ee_{J_\ell},\quad\ee_{J_\ell}\T y^\ell=1,\\
 y^\ell_{h_\ell}<0,\quad
 y^\ell_{h_\ell}\le y^\ell_i\ (i\in J_\ell),\qquad x_h>0.
 \end{gathered}
\]
The chosen support makes complementarity linear; zero entries in the
declared support cause no problem. There are polynomially many variables
and constraints. Replace each strict negative bound by
$y^\ell_{h_\ell}\le-\zeta$, replace $x_h>0$ by $x_h\ge\zeta$, and add
$0\le\zeta\le1$. A real witness gives a feasible point with $\zeta>0$.
Maximizing $\zeta$, the standard rational linear-programming bit bound
provides a positive rational optimum and a witness with polynomial
encoding length. Free variables may be split into nonnegative parts.
The prefix, support, and rational variables are consequently a
polynomial-size certificate; all displayed conditions are checked by
exact rational arithmetic. This proves coNP membership.
\end{proof}

The strict face property is used only for the reduction family. It does
not remove the boundary distinction between full-start and hereditary
correctness established by Theorem~\ref{th:main}.

\appendix
\section{A basic failure of unrestricted minimum-coordinate deletion}
\label{app:basic-failure}
The following example separates feasibility of the terminal affine
solution from optimality. It does not exhibit the matrix-class
separation in Theorem~\ref{th:main}.
Take
\[
 v_1=(1,2,1)\T,\quad v_2=(2,1,1)\T,\quad v_3=(0,4,1)\T,
 \qquad p=0.
\]
These linearly independent vectors form a triangle in the affine plane
$z=1$. The affine projection of $p$ has weights $(8,-4,-3)$, so $\MIN$
first removes $v_2$. The affine weights on $v_1,v_3$ then become
$(8/5,-3/5)$, and the algorithm returns $v_1$, of squared norm $6$.
The actual projection is $(v_1+v_2)/2$, of squared norm $11/2$.
For an exact certificate, let
\[
 A=V\T V=\begin{pmatrix}6&5&9\\5&6&5\\9&5&17\end{pmatrix},
 \qquad x^*=(1/2,1/2,0)\T.
\]
The Gram matrix is positive definite and
$Ax^*=(11/2)\ee+(0,0,3/2)\T$, which is the simplex KKT certificate
for $b=0$. Thus the first deletion already removes a positive optimal
coordinate. This verifies the stated limitation of the literal rule
described in \cite[Section~2.2]{cevikalp}, without making a claim about
when this limitation was first observed.

\section{Why adding a return test does not repair the rule}
\label{app:return}

This appendix concerns a modification of $\MIN$, not the irreversible
procedure classified in the main text. Its purpose is to distinguish a
sound stopping test from a convergent algorithm. The cycling construction
below is derived from an old complementarity example; neither existence
of cycling nor its stability is claimed as an independent contribution.

\subsection{The precise rule and its historical source}
On a nonempty working set $J$, compute the affine minimizer $y^J$ and
its multiplier $\nu_J$, and put $r^J=Ay^J-b-\nu_J\ee$.
Start from $J=P$ and repeat:
\begin{enumerate}
\item If $\min_{i\in J}y_i^J\le0$, delete one index attaining this
minimum and recompute.
\item Otherwise, stop if all excluded residuals are nonnegative.
If some are negative, insert an index attaining
$\min_{i\notin J}r_i^J$ and recompute.
\end{enumerate}
This \emph{minimum-weight return rule} gives priority to primal
infeasibility. Unlike $\MIN$, it deletes zero weights and permits returns.
It does not maintain a separate feasible iterate or use a blocking-ratio
test. Those distinctions matter: Clark's complete method
\cite[Sections~2.1 and~2.5]{clark} and feasibility-preserving affine
correction are different procedures.

Stickney--Watson \cite{stickneywatson} studied cycling of Bard-type
complementarity algorithms. We use the explicit SPD example attributed
to that paper by Gao--G\"artner--Lamperski
\cite[Section~5, p.~12 of arXiv version~2]{gaogartner}:
\begin{equation}\label{eq:sw-old}
 M=\begin{pmatrix}5&-10&2\\-10&41&-6\\2&-6&1\end{pmatrix},
 \qquad q=(1,-7,1)\T.
\end{equation}
The leading principal minors are $5,105,1$. Our starting point is the
reproduction in \cite{gaogartner}. The scaling and simplex embedding
that follow are our deductions from those data, rather than a theorem
quoted from the original source. In particular, no identification of its
historical pivot selection rule with the present rule is assumed.

Let $D=\diag(1,1,30)$, $N=DMD$, and $p=Dq$. For the complementarity
system $w=Nz+p$, a basic set $K$ has
$z_K=-N_{KK}^{-1}p_K$, $z_{K^c}=0$, and $w_K=0$.
Direct rational calculation gives the following strict cycle for
primal-priority minimum selection:
\[
 123\longrightarrow13\longrightarrow1\longrightarrow\varnothing
 \longrightarrow2\longrightarrow23\longrightarrow123.
\]
At the first three states, the basic vectors on their supports are
$(3,-1,-13/30)$, $(1,-1/10)$, and $(-1/5)$, respectively.
At the empty set the residual is $(1,-7,30)$.
At $K=\{2\}$ the two excluded residuals are $-29/41$ and $-30/41$;
at $K=\{2,3\}$ the sole excluded residual is $-3/5$.
Thus each indicated choice is unique. Diagonal scaling is needed here:
the unscaled data do not cycle under this particular primal-priority
selection from the full basic set.

\subsection{An exact cycle on a simplex with four coordinates}
Write $x=(x_1,z)$ with $z=(x_2,x_3,x_4)$ and define the invertible map
$R(x_1,z)=(x_1+\ee\T z,z)$. Set $\beta=1/10$ and
\begin{equation}\label{eq:sw-simplex}
 A_{\rm sw}=10R\T
 \begin{pmatrix}1&\beta p\T\\\beta p&N\end{pmatrix}R
 =\begin{pmatrix}
 10&11&3&40\\11&62&-96&641\\
 3&-96&406&-1767\\40&641&-1767&9070
 \end{pmatrix}.
\end{equation}
Since $p\T N^{-1}p=3$, the Schur complement of $N$ in the middle
matrix is $1-3\beta^2=97/100>0$. Hence $A_{\rm sw}$ is SPD;
its leading principal minors are $10,499,103540,8730000$.
On the affine simplex hyperplane the homogeneous objective becomes
\[
 \tfrac12x\T A_{\rm sw}x=5+p\T z+5z\T Nz.
\]
On each face containing coordinate $1$, its other affine weights are
therefore $\beta$ times the corresponding basic vector of the scaled
complementarity problem. Table~\ref{tab:sw-simplex} verifies that the
anchor weight $x_1$ stays positive. The excluded simplex residuals are
exactly the corresponding complementarity residuals, because
$10\beta=1$.

\begin{table}[htbp]
\centering\small
\begin{tabular}{@{}clll@{}}
\toprule
$J$&$y_J^J$&Excluded residuals&Action\\
\midrule
$1234$&$(253,90,-30,-13)/300$&none&delete $3$\\
$124$&$(91,10,-1)/100$&$r_3=1$&delete $4$\\
$12$&$(51,-1)/50$&$(r_3,r_4)=(-5,18)$&delete $2$\\
$1$&$(1)$&$(r_2,r_3,r_4)=(1,-7,30)$&add $3$\\
$13$&$(403,7)/410$&$(r_2,r_4)=(-29,-30)/41$&add $4$\\
$134$&$(1469,30,1)/1500$&$r_2=-3/5$&add $2$\\
\bottomrule
\end{tabular}
\caption{A full-start six-state cycle obtained from the old data
\eqref{eq:sw-old}. Weights are listed in increasing index order.
The multipliers are, respectively,
$97/10,49/5,499/50,10,4051/410,247/25$.}
\label{tab:sw-simplex}
\end{table}

For $b=0$, the return rule consequently follows
\[
 1234\longrightarrow124\longrightarrow12\longrightarrow1
 \longrightarrow13\longrightarrow134\longrightarrow1234.
\]
The unique optimum, which is never reached, is
\[
 x^*=(166/175,29/1050,1/42,0)\T,
 \qquad \nu^*=5177/525,
 \qquad r^*=(0,0,0,26/7)\T.
\]
These data satisfy the KKT equations directly. All pivot comparisons
are strict, so continuity of finitely many affine solutions and
residuals on the SPD cone preserves the cycle on an open neighborhood.
Adding $1768\ee\ee\T$ makes every matrix entry strictly positive and
leaves all affine minimizers and residuals unchanged. These robustness
and positivity observations are elementary consequences of the example.

\subsection{Termination with at most three coordinates}
For completeness, we give a small-dimensional bound for the precise
rule above. We do not assert historical priority for this bound.
If coordinate $h$ is deleted from a face $J$, its residual at the child
affine minimizer is
\begin{equation}\label{eq:return-deleted}
 r_h=-y_h^J/(Q_J)_{hh}.
\end{equation}
Indeed, changing only $b_h$ by $-y_h^J/(Q_J)_{hh}$ makes its affine
weight zero, while leaving the child problem unchanged. Comparing its
residual for the original and changed linear terms proves the formula.

\begin{proposition}\label{prop:return-three}
For $n\le3$, every full-start execution of the minimum-weight return
rule terminates at the unique optimum within five affine solves, for
every SPD $A$, every $b$, and every resolution of ties.
\end{proposition}
\begin{proof}
A nonnegative full affine solution is optimal. Removing zero weights
preserves it and its full KKT conditions, so this cleanup terminates
within $n$ solves. On an edge, removing a nonpositive affine weight
leaves the minimizing endpoint, with nonnegative removed residual by
\eqref{eq:return-deleted}. This proves the assertion for $n\le2$.

For $n=3$, suppose the full affine solution is infeasible and its first
deleted coordinate is $h$. Its weight is strictly negative. If the
child edge solution is nonnegative, its residual at $h$ is positive by
\eqref{eq:return-deleted}; hence it is optimal, including after zero
cleanup. Otherwise the next deletion, of $i$, reaches the remaining
vertex $k$. In coordinates $z_h=x_h$, $z_i=x_i$,
$x_k=1-z_h-z_i$, write the objective, up to its constant, as
\[
 \ell_hz_h+\ell_i z_i+
 \tfrac12(d_hz_h^2+2cz_hz_i+d_i z_i^2),
 \qquad d_h,d_i>0,\quad \Delta=d_hd_i-c^2>0.
\]
The numbers $\ell_h,\ell_i$ are the two excluded residuals at $k$.
The deletion of $i$ gives $\ell_i\ge0$. If $k$ is not optimal, its
only negative residual is $\ell_h=-L$ with $L>0$, and $h$ is inserted.
The full affine weight previously deleted at $h$ was
\[
 y_h^P=(d_iL+c\ell_i)/\Delta<0.
\]
Thus $c=-a<0$ and $\ell_i>d_iL/a$. On the returned edge $hk$ the
weight at $h$ is $t=L/d_h>0$.
If $t\le1$, the edge solution is feasible and its residual at $i$ is
\[
 \ell_i-aL/d_h>d_iL/a-aL/d_h=\Delta L/(ad_h)>0.
\]
It is therefore optimal, with one final zero deletion if $t=1$.
If $t>1$, the rule deletes $k$ and reaches vertex $h$. Its residual
toward $k$ is $L-d_h>0$; the residual toward $i$ satisfies
\[
 \ell_i+L-a-d_h>(d_i/a+1)L-a-d_h
 >d_id_h/a-a=\Delta/a>0.
\]
This vertex is optimal. The longest chain is a full solve, an edge,
a vertex, a returned edge, and its endpoint: five solves. The argument
uses no strict comparison between the initial minimum and another
weight, so it permits ties. A zero second deletion cannot cause a
return: it would give $\ell_i=0$, contradicting $y_h^P<0$ in the
nonoptimal-vertex case. This also covers the boundary cases.
\end{proof}

The example and proposition locate the smallest number of coordinates
for full-start cycling of this precise return rule at four, or geometric
simplex dimension three. A feasible blocking-ratio step is not this
rule: it chooses a boundary hit on a segment from a current feasible
point toward the new affine minimizer, rather than deleting the most
negative coordinate of that affine minimizer. The example therefore
does not contradict convergence results for that correction.

All entries of Table~\ref{tab:sw-simplex} can be checked by the identities
\[
 (A_{\rm sw})_{JJ}y_J^J=\nu_J\ee_J,\qquad
 \ee_J\T y_J^J=1,\qquad
 r_{J^c}^J=(A_{\rm sw})_{J^cJ}y_J^J-\nu_J\ee_{J^c}.
\]
The displayed weights and residuals give every strict comparison in the
cycle. No claim that the scaling and embedding have not previously
appeared elsewhere is made.

\section*{Scope of the results}
The result concerns the stated irreversible minimum-coordinate rule.
It neither asserts a convergence property of other active-set methods
nor identifies the earliest occurrence of related matrix criteria.
First-step safety is polynomially recognizable, whereas both full-start
and hereditary correctness are coNP-complete. Their equal complexity
does not identify the two matrix classes.

\end{document}